\documentclass[11pt,notitlepage,twoside]{article}

\usepackage{latexsym}
\usepackage{amsmath}
\usepackage{amsfonts}
\usepackage{amssymb}

\newcommand{\be}{\begin{equation}}
\newcommand{\ee}{\end{equation}}
\newcommand{\hf}{\hfill $\diamondsuit$}
\newenvironment{pf}{\noindent{\bf Proof.}\enspace}{
\hfill Q.E.D.}

\usepackage{latexsym}
\newtheorem{thm}{Theorem}[section]
\newtheorem{pro}{Proposition}[section]
\newtheorem{lem}{Lemma}[section]
\newtheorem{rem}{Remark}[section]
\newtheorem{cor}{Corollary}[section]
\newtheorem{df}{Definition}[section]
\newtheorem{ex}{Examples}[section]
\numberwithin{equation}{section}
\DeclareFontFamily{T1}{cmr}{\hyphenchar\font=-1}
\begin{document}
\centerline {\bf \Large A CLASS OF    $n$-TUPLES WITH FINITE }
\vskip 0.2 cm
\centerline {\bf \Large  ASCENT AND PERTURBATIONS }
\vskip 1 cm

{\centerline{\bf{Boulbeba ABDELMOUMEN, Hafedh DAMMAK and Sonia
YENGUI}}}
\vskip 0.5 cm
{\centerline{\small { D\'epartement de Math\'ematiques,
Facult\'e des sciences de Sfax, }}}

{\centerline{\scriptsize {Route de Soukra, Km 3.5, BP 1171, 3000,
Sfax, Tunisie.}}}

\vskip 0.5 cm

{\centerline{\small {e-mail: boulbeba.abbelmoumen@ipeis.rnu.tn}}}
{\centerline{\small {e-mail: sonia.yengui@ipeis.rnu.tn}}}
{\centerline{\small {e-mail: dammak-hafedh@yahoo.fr}}}

\vspace{3mm}
\normalsize
\vspace{3mm} \hspace{.05in}\parbox{4.5in}{{\bf\small Abstract. }{\small  \sloppy{
In this paper, we are interested to  study  the stability of the ascent of a mutually commuting $n$-tuple
 $T$ submitted  to
a perturbation by an $n$-tuple  $S$ which commutes mutually with $T$.  This study  lead us to generalize some well known results
for operators and introduce  a class of $n$-tuples in which we obtain the
stability of the ascent of $T$.
}}}

\section{Introduction}
The
problem of stability of the semi-Fredholm operators under additive
perturbations has been around for many years (see for instance
\cite{A.B, A.B1, M.B,   R.G, M.G, T.K}). It is obvious to ask
if these results are extended for some class of $n$-tuples.\vskip 0.2
cm\noindent The ultimate aim of this paper is to study the
stability of the ascent of an $n$-tuple $T=(T_1,...,T_n)$ of mutually commuting operators,
when $T$ is submitted to a perturbation by an $n$-tuple
$S=(S_1,...,S_n)$ which commutes mutually with
$T$  (see
Theorem \ref{theo2}).    In order to make the problem precise,
we introduce  a new class of $n$-tuples, denoted by $\mathcal{K}^e_n(X)$,
which is equal, for $n=1,$ to the set of essentially Kato
operators. In \cite[Proposition
$2.6$]{W.T}, T. West proves that if $T$ is upper semi-Fredholm operator, then \be
\label{(e_1)} a(T)<\infty \Longleftrightarrow
\overline{N^{\infty}(T)}\cap R^{\infty}(T)=\{0\}.\ee In the present
paper,  we prove that the result (\ref{(e_1)}) remains  true for an $n$-tuple which belong to  a class  denoted by $\mathcal{C}_n(X)$ (see Corollary
\ref{cor21}).
In Theorem \ref{theo1}, we prove that, for all $T\in
\mathcal{K}^e_n(X),$ there exists $\varepsilon>0$ such that, for all
 $S$ be in
$\mathcal{B}_n(X)\cap B(O,\varepsilon)$ and mutually commuting with $T$, we have
$$R^{\infty}(T+S)\cap \overline{N^{\infty}(T+S)}\subset
R^{\infty}(T)\cap \overline{N^{\infty}(T)}.$$\noindent This result
is well known  for the case of essentially Kato operators (see
\cite[section $21$, Theorem $14$]{Muller}). Finally, in Theorem \ref{theo2}, we study  the stability of the ascent of a mutually commuting $n$-tuple
 $T$ subjected to
a perturbation by an $n$-tuple  $S$ which commutes mutually with $T$. More precisely, we prove the following:
  $$T\in \mathcal{G}_{n\alpha}(X),\  a(T)<+\infty\mbox{ and }S\in \mathcal{P}(\mathcal{G}_{n\alpha}(X)),\ \alpha>0  \Longrightarrow a(T+S)<+\infty.$$
\noindent This stability result, generalize the
following well known Rako\`{c}evi\'{c} theorem:\vskip 0.2
cm\noindent For $T,K$ be two bounded linear operators such that $TK=KT$. We have
:
$$T\in \Phi_+(X),\ a(T)<\infty\ \textrm{and}\ K\in
\mathcal{P}(\Phi_+(X))\Longrightarrow a(T+K)<\infty.$$
 \vskip 0.2 cm\noindent

Let us introduce some notations.
For $X$ and $Y$ be two Banach spaces we denote by
  $\mathcal{B}(X,Y)$  the set of all bounded linear
operators from   $X$ to $Y$. We write for short $\mathcal{B}(X)=\mathcal{B}(X,X)$.
  By $X^*$ we denote the dual of $X$. Let $M$ be a subset of a Banach space $X$. The annihilator of $M$
is the closed subspace of $X$ defined by
$$M^\perp:=\{f\in X^*;\ f(x)=0,\ \forall x\in M\},$$ while the
pre-annihilator of a subset $W$ of $X^*$ is the closed subspace of
$X$ defined by $$ ^\bot W:=\{x\in X;\ f(x)=0,\ \forall f\in W\}.$$
In the following lemma, we recall some well known results.
 \begin{lem}\label{l0}\cite{Muller} Let $X$ be a Banach space.\vskip 0.2 cm\noindent
$(i)$ If $M$ is a subspace of $X$, then\  $^\bot
(M^\bot)=\overline{M}.$\vskip 0.2 cm\noindent $(ii)$ If
$\{M_\alpha\}_\alpha$ is any family of subsets of $X$, then\
$(\displaystyle
\bigcup_{\alpha}M_\alpha)^\bot=\bigcap_{\alpha}M_\alpha^\bot.$\vskip
0.2 cm\noindent If $\{L_\alpha\}_\alpha$ is any family of subsets
of $X^*$, then \ $\displaystyle
^\bot(\bigcup_{\alpha}L_\alpha)=\bigcap_{\alpha}$ $^\bot
L_\alpha.$\vskip 0.2 cm\noindent $(iii)$ Let $n\in
\mathbb{N}\backslash
\{0,1\}$ and $M_1,..., M_n$ be closed subspaces of $X$. \\
If, $\forall j\in \{1,...,(n-1)\}, \displaystyle\sum_{i=1}^j M_i$
is closed, then $\displaystyle(\sum_{i=1}^n
M_i)^\perp=\bigcap_{i=1}^n M_i^\perp.$\hf
 \end{lem}
 For an operator $T \in \mathcal{B}(X)$ define  $T^*$ the dual operator of $T.$ We write $N(T)$ for the null space and $R(T)$ for the range of $T$.  The generalized range and generalized kernel of $T$ are defined respectively by
$R^\infty(T)=\displaystyle\bigcap_{k=0}^{+\infty}R(T^k)\mbox{ and }N^\infty(T)=\displaystyle\bigcup_{k=0}^{+\infty}N(T^k).$
 $T$  is said to
be essentially Kato if $R(T)$ is closed and   there exists a
finite-dimensional subspace $F\subset X$ such that
$$N^{\infty}(T)\subset R^{\infty}(T)+F.$$\noindent
The set of upper  semi-Fredholm operators is
defined  by $$\Phi_{+}(X)=\{ T\in \mathcal{B}(X);\ dim
N(T)<\infty \ \mbox{and}\ \mathcal{R}(T) \ \mbox{is closed in}\
X \}.$$   The perturbation class associated with $\Phi_+(X)$ is
defined by  $$\mathcal{P}(\Phi_+(X)):=\{T\in \mathcal{B}(X);\ T+S\in
\Phi_+(X),\ \forall S\in \Phi_+(X)\}.$$
\noindent
\vskip 0.3 cm\noindent
Let $T=(T_1,...,T_n)$ be $n$-tuple. $T$ is said to be of mutually commuting operators on a Banach space $X$ if $\forall(i,j)\in \{1,...,n\}^2$, $T_iT_j=T_jT_i$. We denote
$\mathcal{B}_n(X)$  the set of $n-$tuples  of mutually commuting operators on a Banach space $X$.
$T$ is said to be
of mutually commuting with an $n$-tuple $S=(S_1,...,S_n)$ if $\forall
(i,j)\in \{1,...,n\}^2,\ T_iS_j=S_jT_i.$
Define  $T^*=(T_1^*,...,T_n^*)$ the dual of
$T.$
 We use the standard
multiindex notation. We denote for all
$k\in \mathbb{N}$: $$T^k=(T_1^k,...,T_n^k),\ N(T^k)=\bigcap_{i=1}^n
N(T_i^k)\mbox{ and } R(T^k)=\sum_{i=1}^nR(T_i^k).$$
The generalized range and generalized kernel of the $n$-tuple $T$ are defined respectively by
$R^\infty(T)=\displaystyle\bigcap_{k=0}^{+\infty}R(T^k)\mbox{ and }N^\infty(T)=\displaystyle\bigcup_{k=0}^{+\infty}N(T^k).$
The $n$-tuple $T$  is said to
be essentially Kato if $R(T)$ is closed and   there exists a
finite-dimensional subspace $F\subset N^\infty(T) $ such that
$$N^{\infty}(T)\subset R^{\infty}(T)+F.$$\noindent
  We write $\delta_T :
X\rightarrow X^n$ the bounded operator defined  by
$\delta_T(x)=(T_1x,...,T_nx)\ (x\in X).$  The set of upper semi-Fredholm and lower semi-Fredholm
n-tuples of commuting operators on $X$ are defined respectively by
\begin{align*} \Phi_+^{n}(X)&=\{T=(T_1,...,T_n)\in \mathcal{B}_n(X);\
dim N(T)<\infty\ \textrm{and}\ R(\delta_T)\ \textrm{is
closed}\}\cr \Phi_-^{n}(X)&=\{T=(T_1,...,T_n)\in
\mathcal{B}_n(X);\ codim R(T)<\infty\}.\end{align*} \noindent For
$T$ be in $\mathcal{B}_n(X)$, the ascent of $T$ is defined by
$$a(T)=\min\{k; \ N(T^k)=N(T^{k+1})\}.$$ \noindent If no such $k$
exists, then we set $a(T)=\infty.$\vskip 0.3 cm\noindent
The plan of this paper is as follows. In Section
$2$, we exhibit a large  class of $n$-tuples in which the formula in \cite[Proposition
$2.6$]{W.T} remains valid. In
Section $3$, we give quantitative stability result for generalized range and generalized kernel of  $n$ tuple in the class $\mathcal{K}^e_n(X)$. The  main result of this section is Theorem
 \ref{theo1} . In Section $4$, we study the stability of the ascent of an $n$-tuple (see Theorem \ref{theo2}).

\section{Class $\mathcal{C}_n(X)$ of $n$-tuples}
The ultimate of this section is to extend  Proposition $2.6$ in \cite{W.T} to a class of  $n$-tuples of commuting operators.
\begin{pro}\label{propo 1}

Let $T=(T_1,...,T_n)$  be in $\mathcal{B}_n(X)$.
We  assume that the following conditions
hold:
\vskip 0,3cm
\noindent $(H_1)$ $\dim N(T)$ is finite.
\vskip 0,3cm
\noindent $(H_2)$ There exists $m\in \mathbb{N}^*$ such that,  $\forall k\geq m$, $\forall (i, j) \in\{1,...,n\}^2$, $T_i^k(N(T_i^{k+1}))\subset N(T_j)$.
\vskip3mm
   Then\\

\hfill$ \overline{N^{\infty}(T)}\cap R^{\infty}(T)=\{0\} \mbox{ implies that } T \mbox{ has a finite ascent}.$\hfill\hf

\end{pro}
\begin{pf}
For $n\in \mathbb{N}$ we denote by $V_n:= N(T)\cap R(T^n).$
Since $(V_n)_n$ is a  decreasing serie of vector subspaces of finite dimensions, then there exists $m\in \mathbb{N}$ such that $V_k=V_m$, $\forall k\geq m
$.
Let $x\in N(T^{m+1})$.
Using hypothesis $(H_2)$, we obtain $T_i^m(x)\in N(T_j)$, $\forall (i,j) \in \{1,...,n\}^2$ and  then $T_i^m(x)\in \displaystyle\bigcap_{j=1}^nN(T_j)=N(T)$.
Hence, $$\forall i \in \{1,...,n\},\ T_i^m(x)\in N(T)\cap R(T^m)=V_m\subset \overline{N^{\infty}(T)}\cap R^{\infty}(T)=\{0\}.$$
Thus, $x\in N(T^m)$.
Finally, $N(T^m)=N(T^{m+1})$ and therefore $a(T)$ is finite.
\end{pf}

\begin{ex}

For $i\in \mathbb{N}^*$, we consider   the unilateral backward weighted shift operator $T_i$ defined on $l^r(\mathbb{N},\mathbb{C})$, $r\geq 1$, by:
$$T_i((x_n)_n)=(w_nx_n)_n \mbox{ with } w_n=0 \Longleftrightarrow n=i-1.$$
For all $i\in \{1,...,n\}$, we have $$N(T_i^{k+1})=\{0\}^i\times \mathbb{C}\times \{0\}^{\mathbb{N}}.$$

\noindent Let $(x_n)_{n\in \mathbb{N}}\in N(T_i^{k+1})$, then  $(x_n)_{n\in \mathbb{N}}=(0,...,0,x_i,0,....)$.
Thus, $$T_i^k(N(T_i^{k+1}))=\{0\}\subset N(T_j), \ \forall j\in \{1,...,n\}.$$

Moreover, we have $N^{\infty}(T)=\{0\}$.
Hence, according to Proposition \ref{propo 1}, $a(T)$ is finite.
\end{ex}

\begin{rem}

 Let $T=(T_1,...,T_n)\in \mathcal{B}_n(X) $.
We suppose that, $\forall i\in\{1,...,n\}$,  $a(T_i)$ is finite. We denote by $p=\max\{a(T_i), 1\leq i \leq n\}$.
Then we obtain $$N(T^{p+1})=\displaystyle\bigcap_{i=1}^nN(T_i^{p+1})=\displaystyle\bigcap_{i=1}^nN(T_i^{p})=N(T^p)$$
and we conclude that $a(T)$ is finite.
The converse is false, indeed, if we consider $T=(T_1,T_2)$ such that $T_1$ is invertible and $a(T_2)$ is infinite, then we obtain, $$\forall k\in \mathbb{N},\ N(T^k)=N(T_1^k)\cap N(T_2^k)=\{0\}=N(T^{k+1}).$$ Thus,  $a(T)$ is finite.
\end{rem}
\vskip 0.3 cm In the following result, we study the converse of Proposition \ref{propo 1}.

\begin{pro}\label{propo 2}

Let $T=(T_1,...,T_n)$  be in $\mathcal{B}_n(X)$.
We  assume that the following conditions
hold.
\vskip 0,3cm
\noindent $(H_3)$ There exists $A=(A_1,...,A_n)\in \mathcal{B}_n(X)$ such that $A$ is mutually commuting with $T$ satisfying
 $N(A)=\{0\}$.
\vskip 0,3cm
\noindent $(H_{4})$  There exists $m\in \mathbb{N}$ such that, $\forall p\geq m$, $\delta_A(N(T^p)\cap R(T^p))\subset R(\delta_{T^p})$,
where $$\begin{array}{llll}        \delta_A: & X & \longrightarrow & X^n \\
                                            & x& \longmapsto & (A_1x,...,A_nx)
                                                                                                                       \end{array}
$$

            Then\\

\hfill$ T \mbox{ has finite ascent implies that }  \overline{N^{\infty}(T)}\cap R^{\infty}(T)=\{0\}.$\hfill\hf

\end{pro}
\begin{pf}
Let $p\geq \sup(a(T),m)$ and $V_p:=N(T^p)\cap R(T^p)$.
Let $x\in V_p$, then, by hypothesis $(H_4)$,   $\delta_A(x)\in R(\delta_{T^p})$. Then there exists $z\in X$ such that $A_j(x)=T_j^p(z)$, $\forall j\in \{1,...,n\}$.
Since, $\forall i\in \{1,...,n\}$, $A_iT_i=T_iA_i$, then  we get $A_jT_j^p(x)=T_j^{2p}(z)$, $\forall j\in \{1,...,n\}$.
Moreover, we have $x\in N(T^p)$, then $z\in N(T^{2p})=N(T^p)$.
Hence, $$0= T_j^p(z)=A_j(x), \ \forall j\in \{1,...,n\}.$$
Thus, $x\in N(A)=\{0\}$ and therefore $V_p=\{0\}$.
On the other hand, the fact that $R^{\infty}(T)\subset R(T^p)$ and $\overline{N^{\infty}(T)}\subset N(T^p)$, we obtain
$$ \overline{N^{\infty}}(T)\cap R^{\infty}(T)\subset V_p=\{0\}.$$
Finally, $ \overline{N^{\infty}(T)}\cap R^{\infty}(T)=\{0\}$.
\end{pf}

\begin{df}\label{Def 2}

Let $T=(T_1,...,T_n)$ be in $\mathcal{B}_n(X)$. $T$ is said to be in $\mathcal{C}_n(X)$ if it satisfies the following hypotheses $(\mathcal{H})$:

$$(\mathcal{H})\left\{
    \begin{array}{l}
      (i)\ T\in \Phi_+^n(X); \\
      \\
      (ii)\mbox{ There exists } k_0\in \mathbb{N}\mbox{ such that, } \forall k\geq k_0, \forall (i,j)\in \{1,...,n\}^2, i\neq j, \\ \\  \quad \quad\quad\quad\quad\quad\quad\quad\quad\quad T_i^k(N(T_i^{k+1}))\subset N(T_j);\\
      \\
      (iii) \mbox{ There exists } m\in \mathbb{N}\mbox{ such that, } \forall p\geq m, \mbox{ there exists }  A=(A_1,...,A_n)\in \mathcal{B}_n(X)\\ \mbox{ mutually commuting with  }
       T \mbox{ satisfying: }\\ \\ \quad \quad\quad\quad\quad\quad N(A)=\{0\}  \mbox{ and }\delta_A(N(T^p)\cap R(T^p))\subset R(\delta_{T^p}).
  \end{array}\right.
$$
\end{df}
\begin{rem}
Notice that for $n=1$, $T\in \mathcal{C}_1(X)$, if and only if, $T\in \Phi_+(X)$. Indeed, if $T\in \mathcal{C}_1(X)$ then $T\in \Phi_+(X)$. Conversely, if
$T\in \Phi_+(X)$, then the hypothesis $(\mathcal{H})$ is satisfied if we take  $A=id_X$.
\end{rem}
The following example proves  that $\mathcal{C}_n(X)$ is not empty.
\begin{ex}\label{exam 2}
Let $T=(T_1,...,T_n)$ be in $\mathcal{B}_n(X) $ satisfying the following hypotheses:
\vskip 0.3 cm\noindent
$(i)$ $T_1$ is invertible in $\mathcal{L}(X)$.
\vskip3mm\noindent
$(ii)$ For all $i\in \{2,...,n\}$, $T_i\in \Phi_+(X)$ and $a(T_i)<+\infty$.
\vskip4mm\noindent
Then $T\in \mathcal{C}_n(X)$. Indeed, let $i\in \{2,...,n\},$  since $a(T_i)<+\infty$, then there exits $m\in \mathbb{N}$ such that, $\forall k\geq m$,
$N(T_i^k)=N(T_i^m)$. For $k\geq m$, we have $$T_i^k(N(T_i^{k+1}))=T_i^k(N(T_i^{m}))=T_i^{k-m}(T_i^m(N(T_i^{m})))=\{0\}.$$
Thus, $(\mathcal{H})(ii)$ is satisfied.
Moreover, since $T_1$ is invertible, then, $\forall p\geq 1$, $N(T^p)=\{0\}$  and therefore $(\mathcal{H})(iii)$ is satisfied.
\end{ex}
\noindent
As a consequence of Propositions \ref{propo 1} and \ref{propo 2}, we get the following result which extends  Proposition
$2.6$ in \cite{W.T}.
\begin{cor}\label{cor21}

Let $T=(T_1,...,T_n)$  be in $\mathcal{C}_n(X)$. We have\\

 \hfill$ T \mbox{ has finite ascent  } \Longleftrightarrow  \overline{N^{\infty}}(T)\cap R^{\infty}(T)=\{0\}.$\hfill\hf
\end{cor}
\section{Class $\mathcal{K}_{n}^{e}(X)$ of $n$-tuples}
\begin{df}\label{Def 3}
Let $\alpha> 0$ and $T=(T_1,...,T_n)\in \mathcal{B}_n(X)$. $T$ is said to be essentially kato-stable if there exists  $\varepsilon > \alpha$ such that,
 for all $S=(S_1,...,S_n)\in \mathcal{B}_n(X)\bigcap B(O,\varepsilon)$ and mutually commuting with $T$, the following assertions holds:
\vskip2mm
$(i)$ $T+S$ is essentially kato i.e. there exists a finite dimensional subspace $F_S$ such that $$N^{\infty}(T+S)\subset R^{\infty }(T+S)+F_S.$$
\noindent$(ii)$ $\displaystyle\sum_{i=1}^n(T^*_i+S^*_i)R^{\infty}(T^*)=R^{\infty}(T^*)$.\\
\noindent$(iii)$ $\displaystyle \prod_{i=1}^n(T_i+S_i)\left(R^{\infty}(T+S) \bigcap N^{\infty}(T+S)\right)=R^{\infty}(T+S) \bigcap N^{\infty}(T+S)$.\\ \noindent
$(iv)$ $\forall k\in \mathbb{N}$, $\forall i\in \{1,...,n\}$, $R(T^*_i+S^*_i)^k$ and $R(T^{*k}_i)$ are closed in $X^*$.\\ \noindent
$(v)$ $\forall k\in \mathbb{N}$, $\forall p\in \{1,...,n\}$, $R((T+S)^k)$, $\displaystyle\sum_{i=1}^pN(T_i^k)$ and
$\displaystyle\sum_{i=1}^pN((T_i+S_i)^k)$ are closed in $X$.
\end{df}
\noindent
We denote by $\mathcal{K}_{n\alpha}^{e}(X)$ the set of all essentially kato-stable $n$-tuples and we consider
$$\mathcal{K}_n^e(X):=\displaystyle\bigcup_{\alpha>0}\mathcal{K}_{n\alpha}^{e}(X).$$
\noindent Hence, we have $T\in \mathcal{K}^e(X)$, if and only if, there exists $\varepsilon>0$ such that, $\forall S=(S_1,...,S_n)\in \mathcal{B}_n(X)\bigcap B(O,\varepsilon)$, the assertions
 $(i)$ - $(v)$ of the previous definition  are satisfied.

\begin{rem}

For $n=1$, we have $T\in \mathcal{K}_1^e(X) $ is equivalent to $T$ is essentially Kato operator.
Indeed, it is clear that if $T\in \mathcal{K}_1^e(X)$, then $T$ is essentially Kato.
Conversely, if  $T$ is essentially Kato, then according \cite[ Theorems 5, 7 and 8 section 21]{Muller},
$T\in \mathcal{K}_1^e(X)$, when $\alpha< \displaystyle\lim_{n\rightarrow +\infty}(\gamma(T^n))^{\frac{1}{n}}$.
\end{rem}
\noindent
The following proposition proves that $\mathcal{K}_n^e(X) $ is a non-empty set.

\begin{pro}\label{propo 3}

Let $T=(T_1,...,T_n)$  be in $\mathcal{B}_n(X)$ such that
$T_1$ is invertible in $\mathcal{B}(X)$ and, $\forall i\in \{2,...,n\}$, $T_i\in \phi_+(X)$.
Then \\

\hfill$ T\in \mathcal{K}_n^e(X).$\hfill\hf
\end{pro}
\begin{pf}

Since $T_1$ is invertible, then there exists $\varepsilon_1$ such that, $\forall S_1\in B(O,\varepsilon_1)$, $T_1+S_1$ is invertible in $\mathcal{B}(X)$.
Let $S=(S_1,...,S_n) $ be such that $\|S\|=\displaystyle \sum_{i=1}^n\|S_i\|<\varepsilon_1$.
Then $\|S_1\|<\varepsilon_1$. Hence, $$N^{\infty}(T+S)\subset N^{\infty}(T_1+S_1) \subset\{0\}$$ and therefore
the assertions $(i)$ and $(iii)$ of Definition \ref{Def 3} are satisfied. Moreover, since $T_1^*$ is invertible,
then $R^{\infty}(T^*)=X^*$ and $(T_i^*+S_i^*)(R^{\infty}(T^*))=X^*$. So the assertion $(ii)$ of Definition \ref{Def 3} is satisfied.
Finally, since, $\forall i\in \{1,...,n\}$, $T_i\in \phi_+(X)$, then there exists $\varepsilon\leq\varepsilon_1$ such that, for all  $S_i\in B(O,\varepsilon)$, $(T_i+S_i)\in \phi_+(X)$. So,    for all $k\in \mathbb{N}$ and all  $j\in \{1,...,n\}$,
 $\displaystyle\sum_{i=1}^jN(T_i+S_i)^k$ and $R(T^*+S^*)^k=X^*$
are subsets of finite dimensional. Thus, $(iv)$ and $(v)$ of Definition \ref{Def 3} are satisfied.
\end{pf}
\vskip 0.2 cm\noindent
In the following Proposition and according to Definition \ref{Def 3}, we deduce some useful   properties in the set $\mathcal{K}_{n\alpha}^e(X)$.

\begin{pro}\label{propo 4}

\noindent $(i)$ If $0<\alpha \leq\beta$, then $\mathcal{K}_{n\beta}^e(X)\subset \mathcal{K}_{n\alpha}^e(X)$.
\vskip 0,2cm
\noindent $(ii)$ For $\alpha>0$ and $\beta>0$, we have $\mathcal{K}_{n\alpha}^e(X)\cap \mathcal{K}_{n\beta}^e(X)=\mathcal{K}_{n\max(\alpha,\beta)}^e(X)$.
\vskip 0,2cm
\noindent $(iii)$ If $T$ is in $\mathcal{K}_{n\alpha}^e(X)$, then $\forall \lambda\in \mathbb{C}^*$ such that $|\lambda|>1$, $\lambda T\in \mathcal{K}_{n\alpha}^e(X)$.\hf
\end{pro}
\begin{pf}
\noindent $(i)$ Let $T\in \mathcal{K}_{n\beta}^e(X)$, then there exists $\varepsilon > \beta$ and hence
 $\varepsilon > \alpha$ such that, $\forall S\in \mathcal{B}_n(X)\cap B(O,\varepsilon)$ and mutually
 commuting with $T$, the assertions $(i)- (v)$ of Definition \ref{Def 3} hold.
\vskip3mm
 \noindent $(ii)$ If $\varepsilon>\alpha$ and $\varepsilon >\beta$,  then $\varepsilon >\max(\alpha,\beta)$.
\vskip3mm
 \noindent $(iii)$ Let $\varepsilon>0$ such that $T\in \mathcal{K}_{n\alpha}^e(X)$. Then for $|\lambda|>1$, we have $\displaystyle\frac{1}{\lambda}S\in \mathcal{B}_n(X)\cap \mathcal{B}(O,\varepsilon)$ and hence $(i)-(v)$ of Definition \ref{Def 3} are satisfied. Thus, $\lambda T\in \mathcal{K}_{n\alpha}^e(X)$.
\end{pf}
\vskip 0.2 cm
\noindent The main result of this section is the following.
 \begin{thm}\label{theo1}
 Let $T=(T_1,...,T_n)$  be in $\mathcal{B}_n(X)$.
 Suppose that $T\in \mathcal{K}_n^e(X)$, then there exists $\varepsilon>0$ such that, for all $S\in \mathcal{B}_n(X)\cap \mathcal{B}(O,\varepsilon)$ and
 mutually commuting with $T$, we have\\

 \hfill$R^{\infty}(T+S)\cap\overline{N^{\infty}(T+S)}\subset R^{\infty}(T)\cap\overline{N^{\infty}(T)}.$\hfill\hf
 \end{thm}
\noindent
To prove this Theorem we shall need two Lemmas.
\vskip2mm
\begin{lem}\label{lemme 1}

Let $T=(T_1,...,T_n)$ and $ S=(S_1,...,S_n)$ be in $\mathcal{B}_n(X)$
 mutually commuting with $T$.
We  assume that the following conditions
hold:
\vskip 0,3cm
\noindent $(A_1)$ $\displaystyle\sum_{i=1}^n(T^*_i+S^*_i)R^{\infty}(T^*)=R^{\infty}(T^*)$.
\vskip 0,3cm
\noindent $(A_2)$  $\displaystyle\sum_{i=1}^pN(T_i^k)$ and $\displaystyle\sum_{i=1}^pN((T_i+S_i)^k)$ are closed in $X$, $\forall (p,k)\in \{1,...,n\}\times \mathbb{N}$.
\vskip 0,3cm
          \noindent $(A_3)$ $R(T_i^{*k})$ and $R((T_i^*+S_i^*)^{k})$ are  closed in $X^*$, $\forall (i,k)\in \{1,...,n\}\times \mathbb{N}$.
            \vskip 0,3cm
            Then \\

            \hfill$\overline{N^{\infty}(T+S)}\subset \overline{N^{\infty}(T)}.$\hfill\hf
\end{lem}
\begin{pf}

\begin{itemize}
  \item We claim that \begin{equation}\label{equ1}
  \overline{N^{\infty}(T+S)}={}^{\perp}\Big(R^{\infty}(T^*+S^*)\Big).\end{equation}
\noindent
By Lemma \ref{l0}, we obtain: \begin{align*}\overline{N^{\infty}(T+S)}  &=
 \ ^{\perp}\Big[\Big(\displaystyle\bigcup_{k=0}^{+\infty}(\bigcap_{i=1}^nN((T_i+S_i)^k))\Big)^{\perp}\Big]=
 \ ^{\perp}\Big[\displaystyle\bigcap_{k=0}^{+\infty}(\bigcap_{i=1}^nN((T_i+S_i)^k))^{\perp}\Big]\cr&= \ ^{\perp}\Big[
         \displaystyle\bigcap_{k=0}^{+\infty}(\bigcap_{i=1}^n \ ^{\perp}R((T_i^*+S_i^*)^{k}))^{\perp}\Big].\end{align*}
\noindent
The fact that $R((T_i^*+S_i^*){k}))$ is closed  $\forall (i,k)\in \{1,...,n\}\times \mathbb{N}$ and by   Lemma  \ref{l0}
we infer that:
\begin{align*}
    \overline{N^{\infty}(T+S)} & = \ ^{\perp}\Big[\displaystyle\bigcap_{k=0}^{+\infty}\Big(\sum_{i=1}^n \ ^{\perp}(R((T_i^*+S_i^*)^{k}))^{\perp})\Big)\Big] =\  ^{\perp}\Big[\displaystyle\bigcap_{k=0}^{+\infty}\Big(\sum_{i=1}^n\overline{R((T_i^*+S_i^*)^{k}))}\Big)\Big] \cr&=  \ ^{\perp}\Big[\displaystyle\bigcap_{k=0}^{+\infty}\Big(\sum_{i=1}^nR((T_i^*+S_i^*)^{k}))\Big)\Big].
  \end{align*}
\noindent
Hence, $\overline{N^{\infty}(T+S)}= \ ^{\perp}(R^{\infty}(T^*+S^*))$ and
 our claim is proved.

  \item We claim that \begin{equation}\label{equ2}
  R^{\infty}(T^*)\subset R^{\infty}(T^*+S^*).\end{equation}
\noindent Let $x\in R^{\infty}(T^*)$. Using $(A_1)$, there exist $ x_1,...,x_n\in R^{\infty}(T^*)$ such that $$x=\displaystyle\sum_{i=1}^n(T_i^*+S_i^*)(x_i).$$
By induction on $ p\in \mathbb{N}^* $, we will have
$$x= \displaystyle\sum_{k_1+...+k_n=np}(T_1^*+S_1^*)^{k_1}...(T_n^*+S_n^*)^{k_n}(z_{k_1,...,k_n}),$$ where $z_{k_1,...,k_n}\in R^{\infty}(T^*).$
Since $k_1+...+k_n=np$, then there exists $ i\in \{1,...,n\}$ such that $k_i\geq p$ and hence   $$(T_1^*+S_1^*)^{k_1}...(T_n^*+S_n^*)^{k_n}(z_{k_1,...,k_n})
\in (T_i^*+S_i^*)^{k_i}(R^{\infty}(T^*))\subset R((T_i^*+S_i^*)^p).$$
Thus, $x\in \displaystyle\sum_{i=1}^nR((T_i^*+S_i^*)^p)$.
Hence, $$x\in \displaystyle\bigcap_{p=0}^{+\infty}(\displaystyle\sum_{i=1}^nR((T_i^*+S_i^*)^p))=R^{\infty}(T^*+S^*),$$ which proves
 our claim.
\end{itemize}

Now, using Equation (\ref{equ2}), we obtain ${}^{\perp}R^{\infty}(T^*+S^*)\subset {}^{\perp}R^{\infty}(T^*)$. Then  Equation (\ref{equ1})  implies that\\

\hfill$\overline{N^{\infty}(T+S)}\subset\overline{N^{\infty}(T)}.$\hfill
\end{pf}

\begin{lem}\label{lemme 2}

Let $T=(T_1,...,T_n)$ be in $\mathcal{B}_n(X)$.
We  assume that, there exists  $\varepsilon > 0$ such that,
 for all $S=(S_1,...,S_n)\in \mathcal{B}_n(X)\bigcap B(O,\varepsilon)$ and mutually commuting with $T$, the following conditions
hold:
\vskip 0,3cm
\noindent $(A_4)$ $T+S$ is essentially Kato and $R((T+S)^k)$ is closed in $X$, $\forall k\in \mathbb{N}$.
\vskip 0,3cm
\noindent $(A_{5})$ $g_{T,S}(R^{\infty}(T+S) \bigcap N^{\infty}(T+S))=R^{\infty}(T+S) \bigcap N^{\infty}(T+S)$, where $g_{T,S}$ is the operator defined by: $$g_{T,S}=\displaystyle \prod_{i=1}^n(T_i+S_i).$$

            Then \\

            \hfill$R^{\infty}(T+S) \bigcap \overline{N^{\infty}(T+S)}\subset R^{\infty}(T).$\hfill\hf
\end{lem}
\begin{pf}
Since  $T+S$ is essentially Kato, then there exists a finite dimensional subspace $F$ such that $F\cap R^{\infty}(T+S)=\{0\}$ and $$\overline{N^{\infty}(T+S)}\cap R^{\infty}(T+S)=(\overline{N^{\infty}(T+S)\cap
R^{\infty}(T+S)}+F)\cap R^{\infty}(T+S).$$ \noindent Since,
$R^{\infty}(T+S)$ is a closed subspace, then $$(\overline{N^{\infty}(T+S)\cap R^{\infty}(T+S)}+F_1)\cap
R^{\infty}(T+S)=\overline{N^{\infty}(T+S)\cap R^{\infty}(T+S)}.$$
\noindent Now, to prove Lemma \ref{lemme 2}, it is sufficient to show that, for $ k\in \mathbb{N }$, $$R^{\infty}(T+S) \bigcap N(T+S)^k\subset R^{\infty}(T).$$\noindent
We will do this by induction on $k$. The statement is clear for $k=0$. Let $k\in \mathbb{N}^*$ and assume that the inclusion holds for $k-1$.
Let $x_0\in R^{\infty}(T+S) \bigcap N(T+S)^k $. Then we have $x_0\in R^{\infty}(T+S)\bigcap N^{\infty}(T+S)$.
Using the hypothesis $(A_5)$ we can find an infinite sequence $x_0,x_1,...$ in $R^{\infty}(T+S)\bigcap N^{\infty}(T+S)$ such that $x_{j-1}=g_{T,S}(x_j)$. By Lemma \ref{lemme 1}, $(x_j)_{j\in \mathbb{N}}\in \overline{N^{\infty}(T)}.$
Since $T$ is essentially Kato, then $\dim(\overline{N^{\infty}(T)}/R^{\infty}(T) \bigcap \overline{N(T)})=m$ is finite.
Thus $x_0,...,x_m$ are linearly dependent, i.e., there exists a non-trivial combination
$$x:= \displaystyle \sum_{i=0}^m\alpha_ix_i\in R^{\infty}(T).$$ let $l$ be such that $\alpha_l\neq0$ and $\alpha_j=0$ for $j=l+1,...,m$.
We obtain: $$g_{T,S}^lx=\alpha_lx_0+\displaystyle\sum_{j=0}^{l-1}\alpha_j g_{T,S}^l x_j\in \alpha_lx_0+(R^{\infty}(T+S) \bigcap N(T+S)^{k-1}).$$
\noindent
Finally, we obtain $x_0\in R^{\infty}(T) $.
\end{pf}

\vskip 0.3 cm\noindent
\noindent {\bf Proof  of Theorem \ref{theo1}}
\vskip3mm\noindent
By Lemma \ref{lemme 1} we have $\overline{N^{\infty}(T+S)}\subset \overline{N^{\infty}(T)}$, then, the use of Lemma \ref{lemme 2} leads to
$$R^{\infty}(T+S)\cap\overline{N^{\infty}(T+S)}\subset R^{\infty}(T)$$
and the theorem is proved.
\hfill
Q.E.D.

 \section{Stability of the ascent of $n$-tuple}
\noindent
 In the following we consider some classes of $n$-tuples commuting operators.
\vskip4mm\noindent
 $(i)$ For $\alpha>0$, consider $\mathcal{G}_{n\alpha}(X):=\mathcal{K}_{n\alpha}^e(X)\cap\mathcal{C}_n(X)$.
\vskip3mm\noindent
 $(ii)$ Let $\mathcal{G}_{n}(X):=\displaystyle \bigcup_{\alpha>0}\mathcal{G}_{n\alpha}(X)$.
\vskip3mm\noindent
 $(iii)$ For $\alpha>0$, let $\mathcal{P}(\mathcal{K}_{n\alpha}^e(X))$ the set of all $n$-tuples $T=(T_1,...,T_n)\in \mathcal{B}_n(X)$ such that, for all $S\in \mathcal{K}_{n\alpha}^e(X)$ and mutually commuting with $T$, $T+S\in \mathcal{K}_{n\alpha}^e(X)$.
\vskip3mm
 \noindent$(iv)$ $\mathcal{P}(\mathcal{G}_{n\alpha}(X))$ the set of all $n$-tuple $T$ such that, for all $S\in \mathcal{G}_{n\alpha}(X)$ and mutually commuting with $T$, $T+S\in \mathcal{G}_{n\alpha}(X)$.
\noindent\vskip3mm
 \noindent$(v)$ Let $\mathcal{P}(\mathcal{G}_{n}(X)):=\displaystyle\bigcap_{\alpha>0}\mathcal{P}(\mathcal{G}_{n\alpha}(X))$.

 \begin{rem}
$(i)$ It is clear that,  for $n=1$, $\mathcal{G}_n(X)=\Phi_+(X)$.
\vskip3mm\noindent
$(ii)$ According to Proposition \ref{propo 3}, the set of $n$- tuples defined in example \ref{exam 2} is included in $\mathcal{G}_n(X)$.
\vskip3mm\noindent
$(iii)$ Notice that for $n=1$ and according to \cite[Theorem 9, Section 21]{Muller} $\mathcal{P}(\mathcal{K}_{n\alpha}^e(X))$ contains the set of all compact and quasinilpotent operators.

\end{rem}
\begin{pro}\label{propo 4}

$(i)$ If $0<\alpha\leq \beta$, then $\mathcal{P}(\mathcal{K}_{n\alpha}^e(X))\subset \mathcal{P}(\mathcal{K}_{n\beta}^e(X))$.
\vskip3mm\noindent
$(ii)$ If $T\in \mathcal{P}(\mathcal{K}_{n\alpha}^e(X))$, then for all $\lambda\in ]0,1]$, $\lambda T \in \mathcal{P}(\mathcal{K}_{n\alpha}^e(X))$.\hf
\end{pro}
\begin{pf}
Follows from Proposition \ref{propo 2} $(i)$ and $(iii)$.\hfill
\end{pf}
\begin{thm}\label{theo5}
 Let $T=(T_1,...,T_n)$  be in $\mathcal{K}_{n\alpha}^e(X)$ and let  $S$ be in $\mathcal{P}(\mathcal{K}_{n\alpha}^e(X))$ and mutually commuting with $T$. Then\\

\hfill$R^{\infty}(T+S)\cap\overline{N^{\infty}(T+S)}\subset R^{\infty}(T)\cap\overline{N^{\infty}(T)}.$\hfill\hf

 \end{thm}
\begin{pf}
Let $\mu\in [0,1]$. By Proposition \ref{propo 4}, $T_{\mu} =T+ \mu S\in \mathcal{K}_{n}^e(X)$.
According to Theorem \ref{theo1}  there exists $\varepsilon>0$ such that, $\forall \lambda\in ]\mu,\mu+\frac{\varepsilon}{\|S\|}[$, we have
$$R^{\infty}(T_{\lambda})\cap\overline{N^{\infty}(T_{\lambda})}\subset R^{\infty}(T_{\mu})\cap\overline{N^{\infty}(T_{\mu})}.$$
\noindent If we apply the above procedure with $\mu_0=0< \mu_1 <...<\mu_p=1$, then,  for all $i\in \{0,...,p-1\}$, there exists $\varepsilon_i>0$
such that  $\mu_{i+1}\in ]\mu_i,\mu_i+\frac{\varepsilon_i}{\|S\|}[$ and  we have $$R^{\infty}(T_{\mu_{i+1}})\cap\overline{N^{\infty}(T_{\mu_{i+1}})}\subset R^{\infty}(T_{\mu_i})\cap\overline{N^{\infty}(T_{\mu_i})}.$$Thus,
$$R^{\infty}(T_{\mu_p})\cap\overline{N^{\infty}(T_{\mu_p})}\subset R^{\infty}(T_{\mu_{p-1}})\cap\overline{N^{\infty}(T_{\mu_{p-1}})}\subset...\subset R^{\infty}(T_{\mu_0})\cap\overline{N^{\infty}(T_{\mu_0})},$$\noindent which implies that\\

\hfill$R^{\infty}(T+S)\cap\overline{N^{\infty}(T+S)}\subset R^{\infty}(T)\cap\overline{N^{\infty}(T)}.$\hfill
\end{pf}
\vskip 0.4 cm\noindent
Now, we are ready to state the main result of this paper.

\begin{thm}\label{theo2}
 Let $T=(T_1,...,T_n)\in \mathcal{B}_n(X)$. Then\\

  \hfill$T\in \mathcal{G}_{n\alpha}(X),\  a(T)<+\infty\mbox{ and }S\in \mathcal{P}(\mathcal{G}_{n\alpha}(X)),\ \alpha>0  \Longrightarrow a(T+S)<+\infty.$\hfill\hf

 \end{thm}
\begin{pf}
 Applying Theorem \ref{theo5} , we have $$R^{\infty}(T+S)\cap\overline{N^{\infty}(T+S)}\subset R^{\infty}(T)\cap\overline{N^{\infty}(T)}.$$
According to Corollary \ref{cor21} we infer that $$R^{\infty}(T+S)\cap\overline{N^{\infty}(T+S)}=\{0\}.$$ Thus, by  the same corollary we get $a(T+S)<+\infty$.
\end{pf}

 \begin{cor}
Let $T=(T_1,...,T_n)\in \mathcal{B}_n(X)$. Then\\

  \hfill$T\in \mathcal{G}_{n}(X),\  a(T)<+\infty\mbox{ and }S\in \mathcal{P}(\mathcal{G}_{n}(X))  \Longrightarrow a(T+S)<+\infty.$\hfill\hf
 \end{cor}
\begin{pf}
The result follows since $\mathcal{G}_{n}(X)=\displaystyle\bigcup_{\alpha>0}\mathcal{G}_{n\alpha}(X)$.
\end{pf}
\begin{rem}

Observe that for $n=1$, we obtain the well known Theorem 1 established by V. Rakocevic in [7] indeed we have $\mathcal{G}_1(X)=\Phi_+(X)$.
\end{rem}


\end{document}